\documentclass[11pt]{amsart} 
\usepackage{amscd,amssymb,amsxtra}
\usepackage{mathrsfs}  
\DeclareSymbolFontAlphabet{\mathrsfs}{rsfs}
\usepackage[mathscr]{eucal} 
\usepackage{comment}
\usepackage{color}
\usepackage{enumitem} 
\usepackage[colorlinks,backref=page,citecolor=refkey]{hyperref}
\usepackage{aliascnt}
\usepackage{booktabs}

\usepackage{marginnote}

\makeatletter
\let\@secnumfont\bfseries
\def\section{\@startsection{section}{1}%
  \z@{4\linespacing\@plus\linespacing}{\linespacing}%
  {\bfseries\centering}}
\def\introsection{\@startsection{section}{1}%
  \z@{3\linespacing\@plus\linespacing}{\linespacing}%
  {\bfseries\centering}}
\def\subsection{\@startsection{subsection}{2}%
   \z@{1.25\linespacing\@plus.7\linespacing}{.5\linespacing}%
   {\normalfont\bfseries}}
\def\subsectionsinline{\def\subsection{\@startsection{subsection}{2}%
  \z@{1\linespacing\@plus.7\linespacing}{-.5em}%
  {\normalfont\bfseries}}}

\makeatother

\numberwithin{equation}{section}

\newcommand{\mynewtheorem}[2]{
  \newaliascnt{#1}{equation}
  \newtheorem{#1}[#1]{#2}
  \aliascntresetthe{#1}
  \expandafter\def\csname #1autorefname\endcsname{#2}
}

\theoremstyle{definition}
\mynewtheorem{definition}{Definition}
\mynewtheorem{example}{Example}
\mynewtheorem{problem}{\color{blue}Problem}
\mynewtheorem{probsec}{\color{blue}Problem}
\mynewtheorem{exercise}{Exercise}
\mynewtheorem{question}{\color{blue}Question}
\mynewtheorem{project}{\color{blue}Project}
\mynewtheorem{construction}{Construction}
\mynewtheorem{task}{\color{blue}Task}
\mynewtheorem{otask}{\color{green}Obsolete Task}
\mynewtheorem{notation}{Notation}

\newtheorem*{definition*}{Definition}
\newtheorem*{example*}{Example}
\newtheorem*{problem*}{\color{blue}Problem}
\newtheorem*{probsec*}{\color{blue}Problem}
\newtheorem*{exercise*}{Exercise}
\newtheorem*{question*}{\color{blue}Question}
\newtheorem*{project*}{\color{blue}Project}
\newtheorem*{construction*}{Construction}
\newtheorem*{notation*}{Notation}

\theoremstyle{remark}
\mynewtheorem{note}{Note}
\mynewtheorem{remark}{Remark}
\mynewtheorem{data}{Data}

\newtheorem*{note*}{Note}
\newtheorem*{remark*}{Remark}
\newtheorem*{data*}{Data}

\theoremstyle{plain}
\mynewtheorem{theorem}{Theorem}
\mynewtheorem{corollary}{Corollary}
\mynewtheorem{lemma}{Lemma}
\mynewtheorem{proposition}{Proposition}
\mynewtheorem{conjecture}{Conjecture}
\mynewtheorem{claim}{Claim}
\mynewtheorem{proposal}{Proposal}
\mynewtheorem{conclusion}{Conclusion}
\mynewtheorem{hypothesis}{Hypothesis}
\mynewtheorem{assumption}{Assumption}

\newtheorem*{theorem*}{Theorem}
\newtheorem*{corollary*}{Corollary}
\newtheorem*{lemma*}{Lemma}
\newtheorem*{proposition*}{Proposition}
\newtheorem*{conjecture*}{Conjecture}
\newtheorem*{claim*}{Claim}
\newtheorem*{proposal*}{Proposal}
\newtheorem*{conclusion*}{Conclusion}
\newtheorem*{hypothesis*}{Hypothesis}
\newtheorem*{assumption*}{Assumption}

\newenvironment{proof*}[1][\proofname]{
  \begin{proof}[#1]}{  
\end{proof}}

\definecolor{refkey}{rgb}{0,.6,.4}

\renewcommand{\:}{\colon}
\renewcommand{\AA}{{\mathbb A}}
\newcommand{\Ahat}{{\hat A}}

\newcommand{\CC}{{\mathbb C}}
\newcommand{\CP}{{\mathbb C\mathbb P}}

\newcommand{\EE}{\mathbb E}

\DeclareMathOperator{\End}{End}

\newcommand{\HH}{{\mathbb H}}

\newcommand{\PP}{{\mathbb P}}
\DeclareMathOperator{\pt}{pt}
\newcommand{\QQ}{{\mathbb Q}}
\newcommand{\RP}{{\mathbb R\mathbb P}}
\newcommand{\RR}{{\mathbb R}}
\newcommand{\TT}{\mathbb T}
\DeclareMathOperator{\Spin}{Spin}

\newcommand{\ZZ}{{\mathbb Z}}

\newcommand{\chiup}{\raise.5ex\hbox{$\chi$}}

\DeclareMathOperator{\ind}{ind}
\newcommand{\inv}{^{-1}}
\DeclareRobustCommand{\mstrut}{^{\vphantom{1*\prime y\vee M}}}

\newcommand{\temsquare}{\raise3.5pt\hbox{\boxed{ }}}

\newcommand{\zmod}[1]{\ZZ/#1\ZZ}

\newcommand{\zt}{\zmod2}

\newcommand{\longhookrightarrow}{\lhook\joinrel\longrightarrow}
\newcommand{\hneg}{\mkern-.5\thinmuskip}
 
\DeclareFontFamily{U}{mathx}{}
\DeclareFontShape{U}{mathx}{m}{n}{<-> mathx10}{}
\DeclareSymbolFont{mathx}{U}{mathx}{m}{n}
\DeclareMathAccent{\widehat}{0}{mathx}{"70}
\DeclareMathAccent{\widecheck}{0}{mathx}{"71}
 
\DeclareMathSymbol{\bigtimes}{1}{mathx}{"91}

\DeclareMathOperator{\SO}{SO}
\DeclareMathOperator{\Sp}{Sp}
\let\O\relax
\DeclareMathOperator{\O}{O}

\DeclareMathOperator{\SU}{SU}
\DeclareMathOperator{\GL}{GL}

\newcommand{\upplus}{^{>0}}
\newcommand{\upnn}{^{\ge0}} 
\newcommand{\Zp}{\ZZ\upplus}
\newcommand{\Znn}{\ZZ\upnn}

\usepackage[all,2cell]{xy}
\usepackage{graphicx}\usepackage{epsf}

\newcommand{\bmuu}{\mbox{$\raisebox{-.07em}{\rotatebox{9.9}
  {\tiny {\bf /}
  }}\hspace{-0.49em}\mu\hspace{-0.88em}\raisebox{-0.98ex}{\scalebox{2} 
  {$\color{white}\phantom{.}$}}\hspace{-0.416em}\raisebox{+0.88ex}
  {$\color{white}\phantom{.}$}\hspace{0.46em}$}} 
\newcommand{\bmut}{\bmu 2}
\newcommand{\bmu}[1]{\bmuu _{#1}}

\newcommand{\lact}{\;\reflectbox{\rotatebox[origin=c]{+90}{$\circlearrowleft$}}\;}  
\newcommand{\ract}{\;\rotatebox[origin=c]{+90}{$\circlearrowleft$}\;} 

\newcommand{\rightarrowdbl}{\rightarrow\mathrel{\mkern-14mu}\rightarrow}
\newcommand{\xrightarrowdbl}[2][]{%
  \xrightarrow[#1]{#2}\mathrel{\mkern-14mu}\rightarrow}

\definecolor{refkey}{rgb}{0,.8,.2}\definecolor{labelkey}{rgb}{1,0,0}

\DeclareMathOperator{\Aff}{Aff}
\DeclareMathOperator{\Cliff}{Cliff}
\DeclareMathOperator{\Pin}{Pin}
\DeclareMathOperator{\Tors}{Tors}
\DeclareMathOperator{\proj}{proj}
\newcommand{\Cnn}{\Cliff_{-n}}
\newcommand{\Cno}{\Cliff_{n,1}}

\newcommand{\Cnp}{\Cliff_{+n}}
\newcommand{\Cpq}{\Cliff_{p,q}}
\newcommand{\GLnR}{\GL_n\!\RR}
\newcommand{\Gnln}{(G_n,\lambda _n)}
\newcommand{\HP}{\HH\PP}
\newcommand{\OP}{\Omega ^{\Pp}}
\newcommand{\Pmn}{\Pin^{-}_n}
\newcommand{\Pm}{\Pin^{-}}
\newcommand{\Ppm}{\Pin^{\pm}}
\newcommand{\Ppn}{\Pin^{+}_n}
\newcommand{\Ppq}{\Pin_{p,q}}
\newcommand{\Pp}{\Pin^{+}}
\newcommand{\QZ}{\QQ/\ZZ}
\newcommand{\RZ}{\RR/\ZZ}

\newcommand{\cE}{\widecheck{E}}
\newcommand{\cH}{\widecheck{H}}
\newcommand{\cKO}{\widecheck{KO}}
\newcommand{\cK}{\widecheck{K}}
\newcommand{\cpM}{\cp^M}
\newcommand{\cp}{\widecheck{\pi}_!}
\newcommand{\cy}{\relax} 
 
\newcommand{\hH}{\widecheck{H}}
\newcommand{\hM}{\widehat{M}}
\newcommand{\mc}{\mathfrak{m}_c}

\newcommand{\sB}{\mathscr{B}}
\newcommand{\sE}{\mathscr{E}}
\newcommand{\sG}{\mathscr{G}}
\newcommand{\sx}{\sigma \mstrut _{\hneg\xi }}
\newcommand{\tKO}{\widetilde{KO}}
\newcommand{\tla}{\tilde{\lambda }}
\newcommand{\ts}{\tilde\sigma }
\newcommand{\w}{\relax} 
\newcommand{\xME}{\xi \mstrut _M(E)}
\newcommand{\y}{\relax}

\begin{document}

\abovedisplayskip18pt plus4.5pt minus9pt
\belowdisplayskip \abovedisplayskip
\abovedisplayshortskip0pt plus4.5pt
\belowdisplayshortskip10.5pt plus4.5pt minus6pt
\baselineskip=15 truept
\marginparwidth=55pt

\makeatletter
\renewcommand{\tocsection}[3]{%
  \indentlabel{\@ifempty{#2}{\hskip1.5em}{\ignorespaces#1 #2.\;\;}}#3}
\renewcommand{\tocsubsection}[3]{%
  \indentlabel{\@ifempty{#2}{\hskip 2.5em}{\hskip 2.5em\ignorespaces#1%
    #2.\;\;}}#3} 
\renewcommand{\tocsubsubsection}[3]{%
  \indentlabel{\@ifempty{#2}{\hskip 5.5em}{\hskip 5.5em\ignorespaces#1%
    #2.\;\;}}#3} 

\def\@makefnmark{%
  \leavevmode
  \raise.9ex\hbox{\fontsize\sf@size\z@\normalfont\tiny\@thefnmark}} 
\def\multfoot{\textsuperscript{\tiny\color{red},}}
\def\footref#1{$\textsuperscript{\tiny\ref{#1}}$}

\providecommand*{\xhookrightfill@}{%
  \arrowfill@{\lhook\joinrel\relbar}\relbar\rightarrow
}
\providecommand*{\xhookrightarrow}[2][]{%
  \ext@arrow 0395\xhookrightfill@{#1}{#2}%
}
\makeatother

\setcounter{tocdepth}{2}



 \title[Index theory on Pin manifolds]{Index theory on Pin manifolds} 

 \author[D. S. Freed]{Daniel S.~Freed}
 \address{Harvard University \\ Department of Mathematics \\ Science Center
Room 325 \\ 1 Oxford Street \\ Cambridge, MA 02138}
 \email{dafr@math.harvard.edu}

 \thanks{This material is based upon work supported by the National Science
Foundation under Grant Number DMS-2005286 and by the Simons Foundation Award
888988 as part of the Simons Collaboration on Global Categorical Symmetries.
This work was performed in part at Aspen Center for Physics, which is supported
by National Science Foundation grant PHY-2210452.}

 \dedicatory{To Is}
 \date{July 23, 2024}
 \begin{abstract} 
 We give a systematic treatment of index theory on Pin manifolds, emphasizing
the Clifford linear Dirac operator and differential $KO$-theory.  This
expository article is based on joint work with Mike Hopkins.
 \end{abstract}
\maketitle

Is Singer was a maverick mathematician, always forward-looking and open to
ideas from all parts of mathematics and beyond.  In tribute to his research,
mentorship, and leadership, I chose to tell a story that touches on different
parts of his mathematics.  As befits Is, its origins lie in physics and the
journey weaves through algebra, differential geometry, global analysis, and
homotopy theory; we end it here with prospects for further development. 
 
In two recent works~\cite{FH1,FH2} that apply geometric methods and homotopy
theory to string theory, quantum field theory, and condensed matter physics,
Mike Hopkins and I ran into Atiyah-Patodi-Singer $\eta $-invariants on Pin
manifolds, that is, manifolds which are unoriented yet admit spinor fields
and Dirac operators.  Index theory on Pin manifolds has appeared sporadically
over the years, for example in~\cite{G,St,Z}, but as far as I know there is no
systematic treatment.  I take this opportunity to share some of what we
learned.  Our exposition places differential $KO$-theory at the center, in part
to promote its role in geometric index theory and also in recognition of Is'
contribution~\cite{HS} to the foundations of generalized differential
cohomology.
 
The main observation is~\eqref{eq:7}--\eqref{eq:8}, where embeddings of Pin
groups into Spin groups are constructed.  This essentially reduces the Pin case
to the standard Spin case, once we work in the framework of \emph{Clifford
linear Dirac operators}.
 
One might ask: Why Pin manifolds?  We give a conceptual answer in the first
half of~\S\ref{sec:2}, based on general principles of symmetry and structure in
differential geometry.  In physics, one often encounters Pin manifolds in
theories with time-reversal symmetry, and that is its origin here.
 
Our exposition begins in~\S\ref{sec:1} with a quick resum\'e of some spinor
algebra.  Geometric structures and the geometry of frame bundles is the
background for the Dirac operator on Spin and Pin manifolds, as we tell
in~\S\ref{sec:2}.  The passage from topological index theory~(\S\ref{sec:3}) to
geometric index theory~(\S\ref{sec:5}) is a passage from topological $K$-theory
to differential $K$-theory, and it is the real ($KO$) version of both that is
relevant here.  In~\S\ref{sec:4} we give a quick introduction to differential
cohomology as a warmup to differential $KO$-theory.  \autoref{sec:6} is a brief
digression to explain one of the motivating physics problems: an anomaly in
string theory.  Topological formul\ae\ for $\eta $-invariants on Pin manifolds
are given in~\S\ref{sec:7}, which concludes with some speculation about further
developments in differential cohomology theory that might apply here and
elsewhere. 
 
The author warmly thanks the referee for their careful reading and constructive
suggestions.

   \section{Clifford algebras and Pin groups}\label{sec:1}

We recommend the classic~\cite[Part~I]{ABS}. 
 
Fix $p,q\in \Znn$ and set $n=p+q$.  The \emph{Clifford algebra}\footnote{My
sign convention for~$\Cpq$ (vs.~$\Cliff_{q,p}$) differs from some references.}
$\Cpq$ is the unital associative algebra with generators $e_1,\dots ,e_n$
subject to the relations 
  \begin{equation}\label{eq:1}
     \begin{gathered}
     \begin{aligned} e_1^2&=\cdots =e_p^2 &\!\!\!=+1 \\
      e_{p+1}^2&=\cdots=e_{p+q}^2 & \!\!\!= -1\end{aligned}  \\[3pt]
      e_ie_j=-e_je_i\quad (i\neq j)\end{gathered}
  \end{equation}
For convenience we use the notation ($n\in \Zp$) 
  \begin{equation}\label{eq:2}
     \Cnp=\Cliff_{n,0}\qquad \Cnn=\Cliff_{0,n} 
  \end{equation}
The Clifford algebra is $\zt$-graded:  
  \begin{equation}\label{eq:49}
     \Cliff\mstrut _{p,q} = \Cliff^0_{p,q}\oplus \Cliff^1_{p,q}; 
  \end{equation}
the even subspace~$\Cliff^0_{p,q}$ is an ungraded algebra.  We deploy the
\emph{Koszul sign rule} relentlessly~\cite[\S1.1]{DM}.  The linear subspace
of~$\Cpq$ generated by $e_1,\dots ,e_n$ is identified with~$\RR^{p,q}$, with
underlying vector space~$\RR^n$ and an inner product of signature~$(p,q)$:
$\langle \xi ,\eta \rangle=(\xi \eta +\eta \xi )/2$, $\xi ,\eta \in
\RR^{p,q}\subset \Cpq$.  A vector~$\xi \in \RR^{p,q}$ with $\langle \xi ,\xi
\rangle\neq 0$ acts on~$\RR^{p,q}$ by twisted conjugation in~$\Cpq$:
  \begin{equation}\label{eq:3}
     \sx \:\eta \longmapsto -\xi \eta \xi \inv 
  \end{equation}
An easy computation shows that $\sx$ is reflection in the hyperplane orthogonal
to~$\xi $.  Observe that $\sigma \mstrut _{t\xi }=\sx$ for all~$t\in \RR^{\neq
0}$.  The subgroup of units in~$\Cpq$ generated by vectors $\xi \in
\RR^{p,q}\subset \Cpq$ with $\langle \xi ,\xi \rangle=\pm1$ is the \emph{Pin
group} $\Ppq\subset \Cpq$.  The map $\xi \to \sx$ extends to a double cover
$\Ppq\to \O_{p,q}$ of the orthogonal group.  If $n=p+q\ge 2$, then the Pin
group~$\Pin_{p,q}$ has four components if $p,q\neq 0$ and two components
if~$p=0$ or~$q=0$.  (See~\eqref{eq:50} below for the case~$n=1$.)  Unless
$p,q\le1$, the \emph{Spin group} is the identity component of the Pin group;
the Spin and Pin groups are subgroups of the units in the even Clifford
algebra:
  \begin{equation}\label{eq:4}
     \Spin\mstrut _{p,q}\subset   \Pin\mstrut _{p,q}\subset  \Cpq^0. 
  \end{equation}
The Spin group double covers the identity component of the orthogonal
group~$\O_{p,q}$.  In the definite case, the Spin and Pin groups are
\emph{compact} Lie groups.  We use the notations
  \begin{equation}\label{eq:5}
     \Ppn=\Pin\mstrut_{n,0}\qquad \Pmn=\Pin\mstrut_{0,n}\qquad
     \Spin\mstrut_n=\Spin\mstrut_{0,n}\cong \Spin\mstrut_{n,0} 
  \end{equation}
Let $\bmu k\subset \TT\subset \CC^\times $ denote the group of
$k^{\textnormal{th}}$~roots of unity inside the group~$\TT$ of unit norm
complex numbers.  In low dimensions there are special isomorphisms of Spin and
Pin groups: 
  \begin{equation}\label{eq:50}
     \begin{tabular}{ c@{\hspace{3.75em}} c@{\hspace{3.75em}}
        c@{\hspace{3.75em}} c} \toprule
        n&\multicolumn{3}{l}{\hspace{-1em}$\Spin_n$\hspace{5.5em} $\Pp_n$
        \hspace{7.4em}$\Pm_n$}\\ \midrule\\[-8pt] $1$ & $\bmut$ &
        $\bmut\times \bmut$ & $\bmu4$ \\ [3pt] $2$ &
        $\TT$&$\bmut\ltimes\TT$&$(\bmu4\ltimes\TT)/\bmut$\\ [3pt] $3$ &
        $\SU_2$&$(\bmu4\times \SU_2)/\bmut$&$\bmut\times \SU_2$\\ [3pt]
        \bottomrule \end{tabular} 
  \end{equation}
 
The $\zt$-graded algebras~$A,A'$ are \emph{Morita equivalent} if there exist
$\zt$-graded vector spaces~$V,V'$ and an isomorphism 
  \begin{equation}\label{eq:6}
     A\otimes \End V\cong A'\otimes \End V' 
  \end{equation}
of $\zt$-graded algebras.  If so, the category of left $A$-modules is
equivalent to the category of left $A'$-modules.  The Clifford algebra
$\Cliff_{1,1}\cong \End\RR^{1|1}$ is Morita trivial; it is Morita equivalent to
the ground field~$\RR$.  We also use the isomorphism $\Cliff_{p',q'}\otimes
\Cliff_{p'',q''}\cong \Cliff_{p'+p'',q'+q''}$.  Finally, $\Cliff_{q,p}$ is the
{opposite algebra} to~$\Cpq$.  (Observe that the Koszul sign rule implies that
multiplication~$*$ in the opposite algebra satisfies
$a_1*a_2=(-1)^{|a_1||a_2|}a_2a_1$.  Apply to generators~$e_i$ in~\eqref{eq:1}.)
 
The keys to index theory on Pin manifolds are the following embeddings: 
  \begin{equation}\label{eq:7}
 \begin{aligned} {\Pp_n}&\longhookrightarrow \Spin_{n,1} &&\!\!\!\subset
  {\Cliff_{n,1}^0} &&\!\!\!\cong \bigl[\Cliff_{+n}\,\otimes \,\Cliff_{-1}\bigr]^0\\ 
      {\Pm_n}&\longhookrightarrow \Spin_{n+1}&&\!\!\!\subset {\Cliff_{+(n+1)}^0}&&\!\!\!\cong
\bigl[\Cliff_{+n}\,\otimes\, \Cliff_{+1}\bigr]^0\end{aligned}  
  \end{equation}
Let $\xi \in \RR^n$ be a unit norm vector, regarded as an element
of~$\Cliff_{+n}$.  Then since unit vectors generate the Pin groups, the
embeddings~\eqref{eq:7} are determined by specifying
  \begin{equation}\label{eq:8}
  \xi \longmapsto \xi \otimes e,
  \end{equation}
where $e$~is the generator of $\Cliff_{\mp1}$.  In~\eqref{eq:8} the first~$\xi
$ lies in~$\Ppm_n$ whereas the second~$\xi $ lies in~$\Cnp$.

  \begin{remark}[]\label{thm:1}
 Note the Morita equivalence 
  \begin{equation}\label{eq:9}
     \Cliff_{n,1}\cong \Cliff_{+(n-1)}\otimes \Cliff_{1,1}
     \overset{\textnormal{Morita}}{\simeq } \;\Cliff_{+(n-1)} 
  \end{equation}
Hence $\Pp_n$~embeds in an algebra Morita equivalent to~$\Cliff_{+(n-1)}$ while
$\Pm_n$ embeds in the algebra~$\Cliff_{+(n+1)}$.  Observe the opposite shifts
in~$n$ for the two embeddings~\eqref{eq:7}:
  \begin{equation}\label{eq:10}
     \begin{aligned} \Pp\: &n\longmapsto n-1 \\ \Pm\: &n\longmapsto
      n+1\end{aligned} 
  \end{equation}
As we shall see, these shifts determine the nature of index theory for the two
types of Pin manifolds.
  \end{remark}

There are analogous embeddings $H_n(s)\hookrightarrow \Cnp\otimes D(s)$ for
other compact Lie groups~$H_n(s)$, parametrized by~$s=0,1$ in the ``complex
case'' and by $s=-3,\dots ,4$ in the ``real case'', as indicated in the
following tables. 

  \begin{equation*}
     \begin{tabular}{ c@{\hspace{2em}} l@{\hspace{2em}} c@{\hspace{2em}}
     c@{\hspace{2em}} c@{\hspace{2em}}} 
     \toprule 
     $s$&$\;\;H^c$&$K$&Cartan&$D$\\ \midrule \\[-8pt] 
     $0$&$\Spin^c$&$\TT$&A&$\CC$\\
     $1$&${ \Pin^c}$&$\TT$&{ AIII}&$\Cliff^{\CC}_{-1}$\\
     \bottomrule \end{tabular} 
  \end{equation*}\medskip
  \begin{equation*}
     \begin{tabular}{ c@{\hspace{2em}} l@{\hspace{2em}} c@{\hspace{2em}}
     c@{\hspace{2em}} c@{\hspace{2em}}} 
     \toprule 
     $\phantom{-}s$&$\quad\quad\; H$&$K$&Cartan&$D$\\ \midrule \\[-8pt] 
     $\phantom{-}0$&$\Spin$&$\bmut$&D&$\RR$\\
     ${-}1$&$\Pin^+$&$\bmut$&DIII&$\Cliff_{-1}$\\
     ${-}2$&$\Pin^+\ltimes\mstrut _{\{\pm1\}}\,\TT$&$\TT$&AII&$\Cliff_{-2}$\\
     ${-}3$&$\Pin^-\times\mstrut _{\{\pm1\}}\Sp_1$&$\Sp_1$&CII&$\Cliff_{-3}$\\ 
     $\phantom{-}4$&$\Spin\,\times\mstrut _{\{\pm1\}}\Sp_1$&$\Sp_1$&C&$\HH$\\
     $\phantom{-}3$&$\Pin^+\times\mstrut _{\{\pm1\}}\Sp_1$&$\Sp_1$&CI&$\Cliff_{+3}$\\ 
     $\phantom{-}2$&$\Pin^-\ltimes\mstrut _{\{\pm1\}}\,\TT$&$\TT$&AI&$\Cliff_{+2}$\\
     $\phantom{-}1$&$\Pin^-$&$\bmut$&BDI&$\Cliff_{+1}$\\
     \bottomrule \end{tabular} 
  \end{equation*}
For each line in the table, there is an evident homomorphism $H_n(s)\to \O_n$,
and its kernel ~$K$ is independent of~$n$.  The Cartan label refers to an
associated symmetric space; the association is explained in the references to
the 10-fold way below.  The group $H_n(4)=\Spin_n\times \mstrut _{\pm1}\Sp_1$
is sometimes denoted `$\Spin^h_n$', a quaternionic analog of the
$\Spin^c$~group; see~\cite{L} and the references therein.

The tables are an instance of a ``10-fold way''.  The 10-fold way originated in
a famous paper of Dyson~\cite{D}, and subsequently has had many incarnations,
of which~\cite{D,AZ,HHZ,K,SRFL,FM,KZ,WS} is a small sample.  This particular
incarnation is~\cite[\S9.2.1]{FH1}.  Our subsequent discussion of index theory
has an analog for each of the compact Lie groups~$H_n(s)$.

   \section{Dirac operators}\label{sec:2}

Is learned differential geometry from Chern at the University of Chicago in the
very late 1940s.  Upon his arrival to MIT in~1950, he and Ambrose set about
developing Chern's ideas into a modern treatment.  One core aspect is the
theory of principal bundles and connections, to which Ambrose-Singer
contributed their eponymous theorem~\cite{AmS}.  These ideas were firmly
established in Is' mind when the need for a Dirac operator on Riemannian
manifolds became the route to explain the integrality of the
$\Ahat$-genus~\cite{S}.  The variation we present here---the \emph{Clifford
linear Dirac operator}---is due to Atiyah-Singer: the earliest written account
I know is~\cite[\S II.7]{LM}.  The basic construction is for Spin manifolds,
and a small variation based on the embeddings~\eqref{eq:7} leads to the Dirac
operator on Pin manifolds.  We begin at the beginning: geometric structures
after Klein and Cartan (Chern's teacher), but somewhat generalized.  In part,
we review these foundations to address the question: Why Pin manifolds?

  \begin{definition}[]\label{thm:2}
 Fix $n\in \Znn$.  An $n$-dimensional \emph{linear symmetry type}~$\Gnln$ is a
Lie group~$G_n$ and a homomorphism $\lambda _n\:G_n\to \GLnR$.
  \end{definition}

  \begin{remark}[]\label{thm:3}
 In Cartan's theory of \emph{$G$-structures}~\cite{C} the homomorphism~
$\lambda _n$ is injective.  We make no such requirement; rather, the kernel
$K=\ker\lambda _n$ is the group of \emph{internal symmetries}.
  \end{remark}

  \begin{example}[]\label{thm:4}
 Riemannian geometry comes from the linear symmetry type $\O_n\hookrightarrow
\GLnR$.  Spin geometry comes from the linear symmetry type
$\Spin_n\rightarrowdbl\SO_n\hookrightarrow \GLnR$, which has internal symmetry
group $K=\bmut$.
  \end{example}

A linear symmetry type determines categories of both affine and curved
geometries.  Below we construct Dirac operators associated to both affine and
curved Spin and Pin geometries.  The \emph{affine symmetry type}
$\tla_n\:\sG_n\to \Aff_n$ associated to a linear symmetry type~$\Gnln$ is
defined by pullback:
  \begin{equation}\label{eq:11}
     \begin{gathered} \xymatrix{&&1\ar[d]&1\ar[d]\\ &&{
     K}\ar@{=}[r]\ar[d]&{ K}\ar[d]\\ 1\ar[r] & \RR^n\ar@{=}[d]\ar[r] &
     {\sG_n}\ar@{-->}[r]\ar@{-->}[d]^{\tla_n} & G_n\ar[r]\ar[d]^{\lambda _n}
     &1\\ 1\ar[r] &\RR^n\ar[r]& \Aff_n\ar[r] & \GLnR\ar[r] &1}
     \end{gathered} 
  \end{equation}
The translation group~$\RR^n$ is a normal subgroup of~$\sG_n$.  The
group~$\sG_n$ acts on real affine space~$\AA^n$ by affine symmetries; the
group~$G_n$ acts on the real vector space~$\RR^n$ by linear symmetries.  For
the linear symmetry type $\O_n\hookrightarrow \GLnR$, the group~$\sG_n$ is the
Euclidean group and the corresponding affine geometry is Euclidean geometry.
 
Symmetries in model geometries become structures on smooth manifolds.  The
linear group~$G_n$ is infinitesimal on a manifold, i.e., it encodes first-order
geometry on tangent spaces.  The translation group acts infinitesimally as a
connection.

  \begin{definition}[]\label{thm:5}
 Fix an $n$-dimensional linear symmetry type~$\Gnln$, and suppose $M$~is a
smooth $n$-dimensional manifold.  A \emph{$\Gnln$-structure} on~$M$ is a
pair~$(P,\theta )$ consisting of a principal $G_n$-bundle $P\to M$ and an
isomorphism
  \begin{equation}\label{eq:12}
  \begin{gathered}
     \xymatrix{\sB(M)\ar[rr]^{\y\theta }_{\cong }\ar[dr]_{\GLnR} && \lambda_n
     (P)\ar[dl]^{\GLnR}\\&M} 
  \end{gathered}
  \end{equation}
of principal $\GLnR$-bundles. 
  \end{definition}

\noindent
 Here $\sB(M)\to M$ is the frame bundle and $\lambda _n(P)\to M$ is the
principal $\GLnR$-bundle associated to $P\to M$ via the homomorphism $\lambda
_n\:G_n\to \O_n$.  Usually we want a \emph{differential} $\Gnln$-structure,
which includes a connection~$\Theta $ on $P\to M$ as well.  The associated
connection~$\lambda _n(\Theta )$ on $\lambda _n(P)\to M$ induces a connection
on the frame bundle, so we can define its torsion.  Torsionfree connections
play an important role in the theory.  (For example, their existence is an
integrability condition on a $\Gnln$-structure.)

  \begin{figure}[ht]
  \centering
  \includegraphics[scale=1]{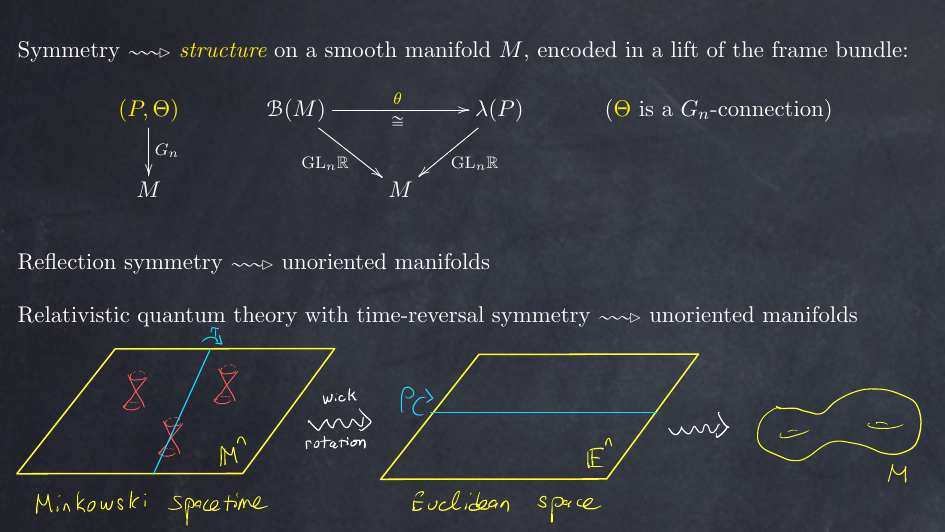}
  \vskip -.5pc
  \caption{\parbox[t]{28pc}{Wick rotation of a spatial reflection in Minkowski
  spacetime to a reflection in Euclidean space; its presence leads to
  unoriented manifolds}}\label{fig:1}
  \end{figure}

Now we can answer the question: Why Pin manifolds?\footnote{Pin manifolds have
associated linear symmetry type $\Pin^{\pm}\rightarrowdbl \O_n\hookrightarrow
\GLnR$.  See~\cite{KT1} for basics on Pin manifolds and for examples in low
dimensions.}  Consider a linear symmetry type~$\Gnln$.  If the image
of~$\lambda _n$ includes orientation-reversing linear maps, then
$\Gnln$-manifolds are unoriented (though they may be orientable); more
precisely, a $\Gnln$-structure does not induce an orientation.  If $\lambda
_n$~factors through the orthogonal group $\O_n\hookrightarrow \GLnR$, then the
affine geometry is Euclidean; if it includes reflection symmetries, then again
$\Gnln$-manifolds are unoriented.  In relativistic quantum field theory, the
symmetry type factors through the group $\O^{\uparrow}_{1,n-1}\hookrightarrow
\GLnR$ of time-orientation preserving isometries of the standard Lorentz
metric.  Relativistic invariance requires that the image of this homomorphism
include the identity component~$\SO^{\uparrow}_{1,n-1}\subset
\O^{\uparrow}_{1,n-1}$.  If the image is the entire
group~$\O^{\uparrow}_{1,n-1}$, then it includes reflections in space.  These
Wick rotate to reflections in Euclidean space (\autoref{fig:1}) and so, as just
mentioned, to unoriented manifolds.  The CPT theorem in quantum field
theory~\cite{SW} implies that such theories have time reflections as well.
Therefore, time-reversal symmetric quantum theories lead to unoriented
manifolds, and in the presence of spinor fields one encounters Pin manifolds.
This is the situation in~\cite{FH1,FH2}.

  \begin{figure}[ht]
  \centering
  \includegraphics[scale=1.3]{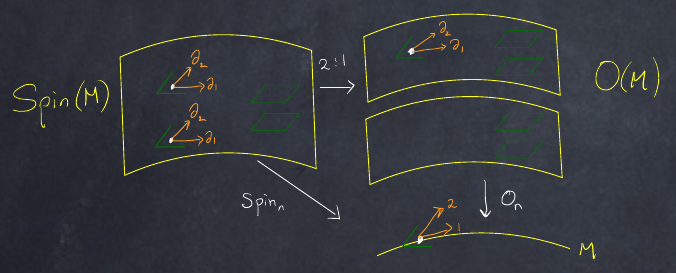}
  \vskip -.5pc
  \caption{\parbox[t]{28pc}{The orthonormal frame bundle, its Levi-Civita
  connection, the tautological horizontal vector fields, and their lift to the
  ``spin frame'' bundle}}\label{fig:2}
  \end{figure}

After these preliminaries we turn to the Dirac operator.  Let $M$~be an
$n$-dimensional Riemannian Spin manifold (linear symmetry type
$\Spin_n\rightarrow\mathrel{\mkern-14mu}\rightarrow\SO_n\hookrightarrow
\GLnR$), and let $\O(M)\to M$ be the principal $\O_n$-bundle of orthonormal
frames.  The fundamental theorem of Riemannian geometry produces a canonical
\emph{Levi-Civita connection} on $\O(M)\to M$.  A point $p\in \O(M)$ in the
fiber at~$m\in M$ is an orthonormal basis $e_1,\dots ,e_n$ of~$T_mM$.  Let
$\partial _1,\dots ,\partial _n$ be the horizontal lifts of the tangent vectors
$e_1,\dots ,e_n$ to~$p$.  The Spin structure is encoded in a principal
$\Spin_n$-bundle $\Spin(M)\to M$, to which the Levi-Civita connection and
horizontal vector fields $\partial _1,\dots ,\partial _n$ lift.  The vector
space\footnote{At this point we raise indices and redefine the algebra $\Cnp$
to have generators $e^1,\dots ,e^n$ with $(e^i)^2=+1$.}  $\Cnp$ carries a left
action of $\Spin_n\subset \Cnp$ that preserves the $\zt$-grading and commutes
with right multiplication $\Cnp\ract\Cnp$.  Let $\gamma ^i$ denote left
Clifford multiplication by~$e^i$.  The Dirac operators on Euclidean
space~$\EE^n$ and on the manifold~$M$ take a similar form:
  \begin{equation}\label{eq:13}
     \begin{aligned} &\EE^n\w:&\y\qquad D & \y=\gamma ^1\frac{\partial
      }{\partial x^1} +\cdots+\gamma ^n\frac{\partial }{\partial x^n}
      &&\lact\quad \w \Bigl(\y\psi \:\EE^n\y\longrightarrow
      \Cnp\w\Bigr)\quad &&\w\ract \quad \Cnp\\[.5pc] &\y M\w:&\y\qquad
      D &\y= \y\,\gamma ^1\,\partial _1 \;\,+\cdots+\,\gamma ^n\,\partial _n
      && \lact\quad \w \Bigl(\y \psi \:\Spin(M)\longrightarrow
      \Cnp\w\Bigr) &&\w\ract \quad \Cnp \end{aligned} 
  \end{equation}
Here a spinor field on~$\EE^n$ is a smooth function $\psi \:\EE^n\to \Cnp$,
while a spinor field on~$M$ is a \emph{$\Spin_n$-equivariant} function $\psi
\:\Spin(M)\to \Cnp$; the equivariance is 
  \begin{equation}\label{eq:14}
     \psi (p\cdot g) = g\inv  \psi (p),\qquad p\in \Spin(M),\quad g\in
     \Spin_n\subset \Cnp. 
  \end{equation}
These Dirac operators are linear maps of right $\Cnp$-modules. 

  \begin{remark}[]\label{thm:6}
 The $\zt$-graded representation $\Spin_n\lact\Cnp$ is reducible.  One can
choose a set of commuting even elements in~$\Cnp$ and simultaneously
diagonalize their action on~$\Cnp$ by right multiplication to pick out an
irreducible subrepresentation; see~\cite[\S8]{AB}.  Instead, we simply record
the entire commuting algebra $\Cnp$.  This is more canonical and more powerful:
it leads naturally to topological and geometric indices of Dirac operators. 
  \end{remark}

  \begin{remark}[]\label{thm:7}
 It is conventional in index theory, certainly since~\cite{ABS}, to work with
\emph{left} $\Cnn$-modules than \emph{right} $\Cnp$-modules.  (Recall that
$\Cnn$ is the opposite $\zt$-graded algebra to $\Cnp$, so these are
equivalent.)  Henceforth, we regard the Dirac operators~\eqref{eq:13} as maps
of left $\Cnn$-modules.
  \end{remark}

The Dirac operator on a $\Pp$ $n$-manifold~$M$ is constructed by a small
modification.  The $\Pp$ structure is encoded in a principal $\Pp_n$-bundle
  \begin{equation}\label{eq:15}
     \Pp_n(M)\longrightarrow \O(M)\longrightarrow M 
  \end{equation}
with lifted Levi-Civita connection and tautological horizontal vector fields
$\partial _1,\dots ,\partial _n$.  The embedding~\eqref{eq:7} produces
commuting actions
  \begin{equation}\label{eq:16}
     \Pp_n\subset \Cno\lact\;\Cno\;\ract\Cno; 
  \end{equation}
the left action of $\Pp_n$ preserves the $\zt$-grading of $\Cno$.  The Dirac
operators are
  \begin{equation}\label{eq:17}
     \begin{aligned} &\EE^n\w:&\y\qquad D & \y=\gamma ^1\frac{\partial
      }{\partial x^1} +\cdots+\gamma ^n\frac{\partial }{\partial x^n}
      &&\lact\quad \w \Bigl(\y\psi \:\EE^n\y\longrightarrow
      \Cno\w\Bigr)\quad &&\w\ract \quad \Cno\\[.5pc] &\y M\w:&\y\qquad
      D &\y= \y\,\gamma ^1\,\partial _1 \;\,+\cdots+\,\gamma ^n\,\partial _n
      && \lact\quad \w \Bigl(\y \psi \:\Pp(M)\longrightarrow \Cno\w\Bigr)
      &&\w\ract \quad \Cno \end{aligned} 
  \end{equation}
These linear maps respect the right $\Cno$-module structures.  Equivalently,
they respect left $\Cliff_{1,n}$-module structures, and by Morita equivalence
we construct a Dirac operator on left $\Cliff_{-(n-1)}$-modules.  It is this
last Dirac operator that we use in the sequel.
 
The case of $\Pm$ manifolds is similar, only now one obtains Dirac operators on
left $\Cliff_{-(n+1)}$-modules. 

  \begin{remark}[]\label{thm:8}
 Note the opposite shifts in~$n$ for $\Pp$ and $\Pm$, already flagged in
\autoref{thm:1}.  Namely, Dirac operator on an $n$-dimensional $\Pp$ manifold
is a map of $\Cliff_{-(n-1)}$-modules, whereas on an $n$-dimensional $\Pm$
manifold the Dirac operator is a map of $\Cliff_{-(n+1)}$-modules.
  \end{remark}

   \section{Topological index theory}\label{sec:3}

We first recall the Atiyah-Singer index theorem for families of Dirac
operators.  The index theorem for general families of elliptic
pseudodifferential operators is proved in~\cite{AS1}; the index of Clifford
linear Fredholm operators is the subject of~\cite{AS2}.  Both papers showcase
Is' deep knowledge and facility with analysis, though that aspect does not come
through in the present account.

  \begin{figure}[ht]
  \centering
  \includegraphics[scale=1.5]{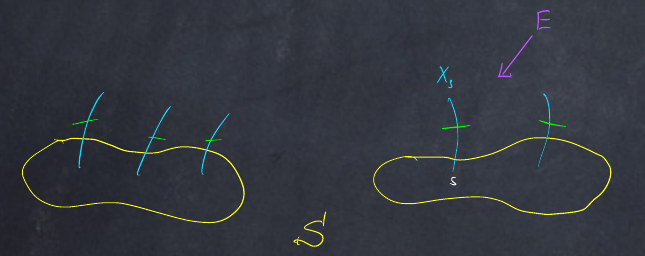}
  \vskip -.5pc
  \caption{A geometric family of Spin manifolds}\label{fig:3}
  \end{figure}

\begin{figure}[ht]
  \centering
  \includegraphics[scale=1.5]{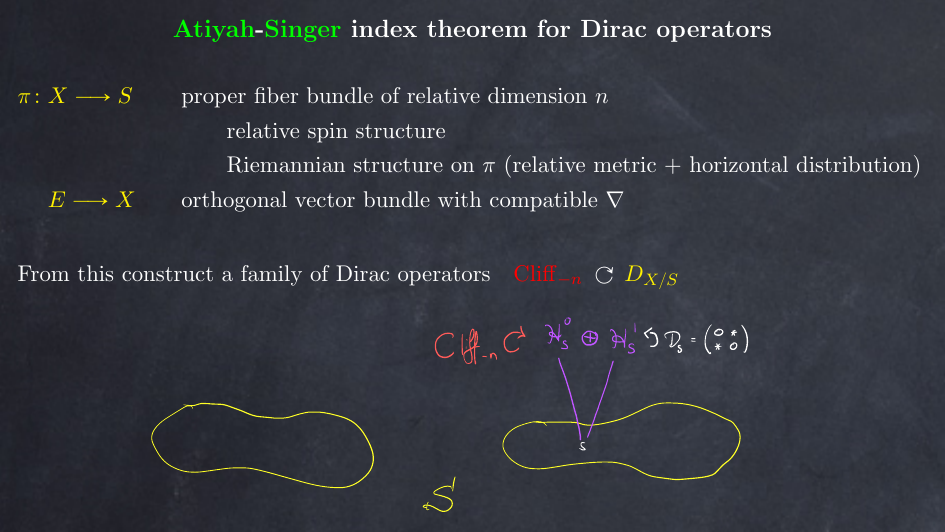}
  \vskip -.5pc
  \caption{The induced family of Clifford linear Dirac operators}\label{fig:4}
  \end{figure}

Let $\pi \:X\to S$ be a proper fiber bundle of smooth manifolds of relative
dimension~$n$, i.e., $S$~smoothly parametrizes a family of closed
$n$-dimensional manifolds.  A \emph{relative Riemannian structure} on~$\pi $ is
a metric on the relative tangent bundle $T(X/S)\to X$ together with a
horizontal distribution on the fiber bundle~$\pi $.  A \emph{relative Spin
structure} is a Spin structure on $T(X/S)\to X$.  Assume given both.  In
addition, suppose $E\to X$ is an orthogonal vector bundle with compatible
covariant derivative.  This data is depicted in \autoref{fig:3}.  From it we
construct a family of Clifford linear Dirac operators~$D_{X/S}$; each $D_s$,
$s\in S$, is a Fredholm map of left $\Cnn$-modules (\autoref{fig:4}).  The
Fredholm, or analytic, index of this family is a class in $KO$-theory:
  \begin{equation}\label{eq:18}
     \ind D_{X/S}\quad \in KO^{-n}(S) 
  \end{equation}
The Atiyah-Singer theorem is a topological formula for this analytic index. 

To push forward along the fibers of $\pi \:X\to S$ in $KO$-theory, at least
over compact subsets of~$S$, choose an embedding $i\:X\hookrightarrow S^N$ into
a sphere and factor the fiber bundle to an embedding followed by projection off
a sphere:
  \begin{equation}\label{eq:19}
     \pi \:X \xhookrightarrow{\pi \times i}S\times
     S^N\xrightarrowdbl{\;\;\proj_1\;\;}S 
  \end{equation}
The key ingredient in the direct image for an embedding, such as $\pi \times
i$, is the construction of the $KO$-theory Thom class~\cite{ABS} of the normal
bundle to the image.  The direct image
  \begin{equation}\label{eq:20}
     \pi _!\:KO^0(X)\longrightarrow KO^{-n}(S) 
  \end{equation}
is the direct image of the embedding followed by the suspension isomorphism in
$KO$-theory.   
 
With this in place, we can state the index theorem. 

  \begin{theorem}[Atiyah-Singer]\label{thm:9}
 $\ind D_{X/S}=\pi _!\bigl([E] \bigr)$, where $[E]\in KO^0(X)$ is the
$KO$-theory class of the real vector bundle $E\to X$. 
  \end{theorem}

The modification to \autoref{thm:9} for families of Pin manifolds is minor; the
main point is the shift~\eqref{eq:10}.  For convenience we restrict to $\Pp$;
the case of $\Pm$ manifolds has the opposite shift but otherwise is the same.
Suppose now that $\pi \:X\to S$ has a relative $\Pp$ structure in place of a
relative $\Spin$ structure.  Then, as explained in~\eqref{eq:17} and
\autoref{thm:8}, there results a family of $\Cliff_{-(n-1)}$-linear Dirac
operators.  Hence the analytic index lies in a shifted $KO$-group:
  \begin{equation}\label{eq:21}
     \ind D_{X/S}\quad \in KO^{-(n-1)}(S) 
  \end{equation}
A similar shift in degree of the Thom class leads to the topological
pushforward map
  \begin{equation}\label{eq:22}
     \pi _!\:KO^0(X)\longrightarrow KO^{-(n-1)}(S) 
  \end{equation}

  \begin{theorem}[]\label{thm:10}
  $\ind D_{X/S}=\pi _!\bigl([E] \bigr)$. 
  \end{theorem}

We illustrate \autoref{thm:10} by noting that the index is a $\Pp$ bordism
invariant: if $M$~is a closed $n$-dimensional $\Pp$ manifold and $\pi ^M\:M\to
\pt$ is the unique map, then $\pi ^M_!(1)\in KO^{-(n-1)}(\pt)$ depends only on
the $\Pp$ bordism class of~$M$ in~$\OP_n$.  The group in which the index takes
its value is determined by the homotopy groups of the orthogonal group, i.e.,
by Bott periodicity:
  \begin{equation}\label{eq:23}
     KO^{-m}(\pt)\;\cong\; \pi _{m-1}O_\infty \;\cong\; \begin{cases} \ZZ,&m\equiv
     0\pmod8;\\ \zt,&m\equiv 1\pmod8;\\ \zt,&m\equiv 2\pmod8;\\ 0,&m\equiv
     3\pmod8;\\ \ZZ,&m\equiv 4\pmod8;\\ 0,&m\equiv 5\pmod8;\\ 0,&m\equiv
     6\pmod8;\\ 0,&m\equiv7\pmod8.\end{cases} 
  \end{equation}

  \begin{example}[$n=2$]\label{thm:11}
 The index of the $\Pp$ Dirac operator is an isomorphism 
  \begin{equation}\label{eq:24}
     \OP_2\xrightarrow{\;\;\cong \;\;}\zt 
  \end{equation}
A generator of $\OP_2$ is the Klein bottle with a nonbounding $\Pp$ structure.
(Note that the real projective plane~$\RP^2$ does not admit a $\Pp$ structure.)
  \end{example}

  \begin{example}[$n=4$]\label{thm:12}
 The bordism group $\OP_4\cong \zmod{16}$ is generated by~$\RP^4$ with either
of its two inequivalent $\Pp$ structures.  However, $KO^{-3}(\pt)=0$ and so the
index provides no information. 
  \end{example}

In the early 1970s index theory pivoted from \emph{global topological}
invariants of general elliptic pseudodifferential operators to \emph{local
geometric} invariants of Dirac operators.  Is was in the vanguard of these
developments: the Atiyah-Patodi-Singer $\eta $-invariant~\cite{APS} is the
first of these geometric invariants.  From a modern perspective, the geometric
invariants are located in \emph{differential} $KO$-theory, to which we now
turn.

   \section{Introduction to differential cohomology}\label{sec:4}

The precursors to differential cohomology are Deligne cohomology~\cite{De} in
algebraic geometry on the one hand and Cheeger-Simons differential
characters~\cite{ChS} on the other, both introduced in the early 1970s.
Together they lead to a smooth version of Deligne cohomology~\cite[\S6.3]{DF},
computed by a cochain complex on a smooth manifold.  Just prior to the turn of
the millennium, different influences from physics\footnote{One motivation is to
encode Dirac charge quantization in the theory of abelian gauge fields.  We
encounter this in~\S\ref{sec:6} in the form of the $C$-field in M-theory.} led
to the need for differential refinements of generalized cohomology
theories~\cite{F1,FH3}.  Characteristically, Is was engaged with these
developments; his joint paper with Hopkins~\cite{HS} lays foundations for the
general theory and has served as the basis for further developments and
applications.  The modern approach, for example in~\cite{BNV}, leans more
heavily on homotopy theory.  The book~\cite{ADH} and recent survey~\cite{Deb}
contain detailed accounts and references.  For differential $K$-theory there
are various geometric models---a sample is~\cite{BS,SS1}---as well as a
refinement of the Atiyah-Singer index \autoref{thm:9} in differential complex
$K$-theory~\cite{FL}.  Recent work on differential $KO$-theory
includes~\cite{GS,GY}.  Our limited goal in this section is to give an entr\'ee
to generalized differential cohomology theory.

We begin with the differential refinement\footnote{Often `$\widehat{H}$' is
used in place of `$\cH$'.  To make explicit the coefficients, one might use
`$\widecheck{H\ZZ}$'.}~$\cH^{\bullet }$ of integer cohomology in low degrees.
Let $M$~be a smooth manifold.  Then
  \begin{equation}\label{eq:25}
     \begin{aligned} {\y H^1(M;\ZZ)} &\cong \bigl\{ \textnormal{smooth maps
      $M\longrightarrow\RZ$} \bigr\} \bigm/ \textnormal{homotopy} \\
      {\y\hH^1(M)}&\cong \bigl\{ \textnormal{smooth maps
      $M\longrightarrow\RZ$} \bigr\} \end{aligned} 
  \end{equation}
The forgetful map $\cH^1(M)\to H^1(M;\ZZ)$ remembers only the homotopy class.
Moving up a degree,
  \begin{equation}\label{eq:26}
     \begin{aligned} {\y H^2(M;\ZZ)} &\cong \bigl\{ \textnormal{principal
     $\RZ$-bundles $P\longrightarrow M$} \bigr\}  \bigm/
     \textnormal{isomorphism} \\ {\y\hH^2(M)}&\cong \bigl\{
     \textnormal{principal $\RZ$-connections $(P,\Theta )\longrightarrow M$}
     \bigr\} \bigm/ \textnormal{isomorphism}  \end{aligned} 
  \end{equation}
The map $\cH^2(M)\to H^2(M;\ZZ)$ forgets the differential geometric data
(connection).  There are various geometric models of ``gerbes'' for degree
three integer cohomology, as well as for connections on gerbes for degree three
differential cohomology.  Such constructions become more intricate as the
degree increases.  For this and many other reasons, one turns to a general
theory.  
 
For any $q\in \Znn$ there is a commutative square of abelian groups 
  \begin{equation}\label{eq:27}
     \begin{gathered}
     \xymatrix@C+1.4pc{\y\hH^q(M)\ar[r]^<<<<<<<<{\textnormal{curvature}}
     \ar[d]_{\pi 
     _0} & \Omega ^q(M)_{\textnormal{closed}}\ar[d]^{\textnormal{de Rham}} \\
     H^q(M;\ZZ)\ar[r]^{} & H^q(M;\RR) } \end{gathered} 
  \end{equation}
The term `curvature' literally applies for~$q=2$ and is used more generally.
Notably, \eqref{eq:27}~is not a \emph{pullback} square; try~$q=1$ and~$M=\pt$.
Rather, $\cH^q(M)$~is constructed as a \emph{homotopy} pullback.  In cochain
complexes, a class in~$\cH^q(M)$ is represented by a triple $(c,h,\omega  )$ in
which  
  \begin{equation}\label{eq:28}
     c\in C^q(M;\ZZ)\qquad h\in
     C^{q-1}(M;\RR) \qquad \omega \in \Omega ^q(M)
  \end{equation}
satisfy 
  \begin{equation}\label{eq:29}
     \delta c=0\qquad\qquad  \delta h=\omega -c \qquad\qquad  d\omega =0
  \end{equation}
The cochain~$h$ serves as a homotopy between~$\omega $ and~$c$, both regarded
as real cocycles of degree~$q$.  \emph{Differential cocycles} $(c,h,\omega )$
are objects of higher groupoids that are \emph{local} on~$M$; see~\cite[\S
A.3]{FN}.

  \begin{figure}[ht]
  \centering
  \includegraphics[scale=1]{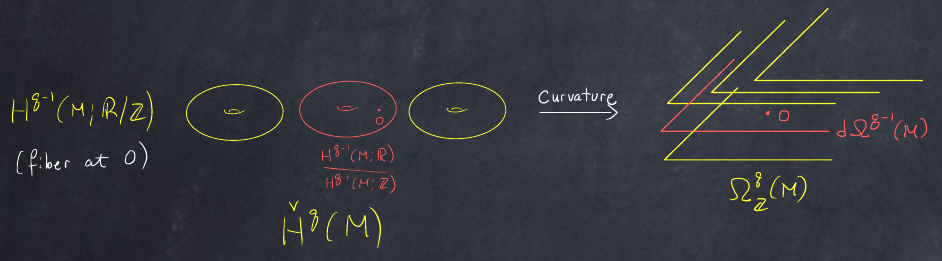}
  \vskip -.5pc
  \caption{The abelian Lie group $\cH^q(M)$}\label{fig:5}
  \end{figure}

There is an abelian Lie group structure on~$\cH^q(M)$, constructed
in~\cite{BSS} and depicted in \autoref{fig:5}.  It reveals many features of
differential cohomology.  The image of the curvature map $\cH^q(M)\to \Omega
^q(M)$ is a union of affine translates of the subspace $d\Omega ^{q-1}(M)$ of
exact $q$-forms.  The set of these translates forms an abelian group isomorphic
to $H^q(M;\ZZ)/\textnormal{torsion}$.  The kernel of the curvature map is the
subgroup of \emph{flat} elements of $\cH^q(M)$.  It is isomorphic to
$H^{q-1}(M;\RZ)$, which is an abelian Lie group whose identity component is the
torus $H^{q-1}(M;\RR)/H^{q-1}(M;\ZZ)$ and whose group of components is
isomorphic to the torsion subgroup $\Tors H^q(M;\ZZ)$.  The group of components
of $\cH^q(M)$ is isomorphic to $H^q(M;\ZZ)$, and the map~$\pi _0$
in~\eqref{eq:27} sends an element of~$\cH^q(M)$ to its component.

  \begin{remark}[]\label{thm:13}
 \ 
 \begin{enumerate}[label=\textnormal{(\arabic*)}]

 \item Topological cohomology features in two homomorphisms involving
differential cohomology:
  \begin{equation}\label{eq:30}
     \begin{alignedat}{2} &{\y \pi _0\:\cH^q(M)\xrightarrowdbl{\;\;\;\;} H^q(M;\ZZ)}
     &\qquad\quad &\textnormal{(underlying topological class)} \\ &{\y \iota
     \:H^{q-1}(M;\RZ)\longhookrightarrow \cH^q(M)} &&\textnormal{(flat
     subgroup)} \end{alignedat} 
  \end{equation}
The first is a quotient and the second is a sub.  The interplay between~$\pi
_0$ and~$\iota $ gives topological information beyond cohomology;
see~\eqref{eq:40} and~\eqref{eq:43} below for examples.

 \item Differential cohomology marries local differential forms with global
integrality. 

 \item There is a calculus of differential cocycles---including products and
integration---that simultaneously lifts (i)~the calculus of differential forms
and (ii)~cup product and pushforward in cohomology.

 \end{enumerate}
  \end{remark}

   \section{Differential $K$-theory and $\eta $-invariants}\label{sec:5}

At some point in the early 1980s, Is suggested three thesis problems to me.
One asked to construct Cheeger-Simons characters for $K$-theory.  Ever
prescient, Is was clearly envisioning differential $K$-theory long before
others.
 
The general theory---as first developed by
Hopkins-Singer~\cite{HS}---constructs $\cK$~and $\cKO$~groups that fit into
commutative squares analogous to~\eqref{eq:27}:
  \begin{equation}\label{eq:31}
     \begin{gathered}
     \xymatrix@R+.5pc@C+1.4pc{\y\cK^q(M)\ar[r]^<<<<<<<<{ 
     \textnormal{\;curvature\;\;\;}}
     \ar[d]_{\pi _0} & \Omega \bigl(M;\RR[u,u\inv ]
     \bigr)^q_{\textnormal{closed}}\ar[d]^{\textnormal{de Rham}} \\
     K^q(M)\ar[r]^{} & H\bigl(M;\RR[u,u\inv ]\bigr)^q } \end{gathered}
     \qquad\quad  (\deg u=2) 
  \end{equation}
and 
  \begin{equation}\label{eq:32}
     \begin{gathered}
     \xymatrix@R+.5pc@C+1.4pc{\y\cKO^q(M)\ar[r]^<<<<<<<<{
     \textnormal{\;curvature\;\;\;}}\ar[d]_{\pi _0} & \Omega
     \bigl(M;\RR[v,v\inv ]\bigr)^q_{\textnormal{closed}}\ar[d]^{\textnormal{de Rham}} \\
     KO^q(M)\ar[r]^{} & H\bigl(M;\RR[v,v\inv ]\bigr)^q } \end{gathered}
     \qquad\quad  (\deg v=4) 
  \end{equation}
In these diagrams $\RR[u,u\inv ]\cong K^{\bullet }(\pt)\otimes \RR$ and
$\RR[v,v\inv ]\cong KO^{\bullet }(\pt)\otimes \RR$.   
 
A useful geometric model represents a differential class by a triple $(E,\nabla
,\eta )$ in which $E\to M$ is a unitary/orthogonal vector bundle with
compatible covariant derivative, and $\eta $~is a differential form of total
degree~$q-1$; see~\cite[\S2]{FL} for details.  Simons-Sullivan~\cite{SS2}
prove, at least for complex $K$-theory, that each differential $K$-theory class
has a representative with~$\eta =0$. 
 
For complex $K$-theory we have isomorphisms 
  \begin{equation}\label{eq:33}
     \begin{aligned} \y&\cK^{-m}(\pt)\xrightarrow[\cong ]{\;\;\pi _0
       \;\;}K^{-m}(\pt)\cong \cy\ZZ,\quad &&\textnormal{$m$~even} \\[6pt]
     \cy\RZ\w\cong 
       K^{-(m+1)}(\pt;\RZ)\xrightarrow[\cong ]{\;\;\iota \;\;}
       \y&\cK^{-m}(\pt),\quad &&\textnormal{$m$~odd}\end{aligned} 
  \end{equation}

  \begin{theorem}[Klonoff~\cite{Kl}]\label{thm:14}
 Let $\y M$~be a closed $n$-dimensional Riemannian Spin${}^c$ manifold, and
suppose given $\y E\to M$ a unitary vector bundle with covariant derivative.
Then the pushforward 
  \begin{equation}\label{eq:35}
     \y\cpM\:\cK^0(M)\longrightarrow \cK^{-n}(\pt) 
  \end{equation}
is computed by the analytic quantities
  \begin{equation}\label{eq:34}
     \y\cpM\bigl([\cE] \bigr)=\begin{cases} \ind D\mstrut _M(E)
     ,&\textnormal{$n$ even}\\\xME\pmod1,&\textnormal{$n$ odd}\end{cases} 
  \end{equation} 
  \end{theorem}

\noindent
 The definition of~\eqref{eq:35} uses the Riemannian structure as well as the
Spin${}^c$ structure; by contrast the direct image in topological $K$-theory
uses only the topological Spin${}^c$ structure.  For $n$~even, the pushforward
is the index of the Dirac operator coupled to the vector bundle $E\to M$,
and~\eqref{eq:34} is the Atiyah-Singer index theorem.  For $n$~odd, we obtain a
geometric invariant---the Atiyah-Patodi-Singer $\eta $-invariant of Dirac
coupled to $E\to M$---but in the form $\xi =\frac{\eta +\dim\ker}{2}\pmod1$, as
in the Atiyah-Patodi-Singer work.\footnote{In the spirit of this paper, we
should use an $\xi $-invariant constructed from the Clifford linear Dirac
operator, as in \cite[\S A.1]{FH2}.  Observe that it is $\xi \pmod1\in \RZ$
that occurs, not $\xi \in \RR$.}  For $n$~odd there is a dependence on metrics
and connections, whereas for $n$~even the result factors though the topological
pushforward.  \autoref{thm:14} is also proved in~\cite{FL} as a corollary to a
differential index theorem for families.

  \begin{remark}[]\label{thm:15}
 \ 
 \begin{enumerate}[label=\textnormal{(\arabic*)}]

 \item Other invariants of geometric index theory---Pfaffian and determinant
line bundles, index gerbes, Bismut superconnection, \dots ---are unified in
differential $K$-theory.

 \item Below we use the extension of \autoref{thm:14} to $\cKO$, though proofs
may not exist in print.

 \end{enumerate}
  \end{remark}

Turning now to differential $KO$-theory, the home for the $\cKO$~invariant of
an $n$-manifold is 
  \begin{equation}\label{eq:36}
     \y\cKO^{-m}(\pt)\w\cong \begin{cases} \RZ ,&m\equiv 7\pmod 8\\
     0,&m\equiv 6\pmod8 \\ 0,&m\equiv 5\pmod8 \\ \ZZ,&m\equiv 4\pmod8 \\
     \RZ,&m\equiv 3\pmod8 \\ \zt,&m\equiv 2\pmod8 \\ \zt,&m\equiv 1\pmod8 \\
     \ZZ,&m\equiv 0\pmod8 \\ \end{cases} 
  \end{equation}
An analog of
\autoref{thm:14} computes the pushforward 
  \begin{equation}\label{eq:37}
     \y\cpM\:\cKO^0(M)\longrightarrow \cKO^{-n}(\pt) 
  \end{equation}
for $M$~is a closed Riemannian Spin manifold.  For $n\equiv 0,1,2,4\pmod8$ we
obtain the index invariants of Atiyah-Singer, including mod~2 indices.  For
$n\equiv 3,7\pmod8$ we obtain $\eta $-invariants.  
 
Now we arrive at the main point, namely the variation---with a
shift~\eqref{eq:10}---for $\Pp$ manifolds.  For safety, we leave it as a
conjecture since---to our knowledge---no proof exists in the literature.

  \begin{conjecture}[]\label{thm:16}
 Let $\y M$~be a closed $n$-dimensional $\Pp$ manifold, and suppose $\y E\to M$
is an orthogonal vector bundle with covariant derivative.  Then 
  \begin{equation}\label{eq:38}
     \y\cpM\bigl([\cE] \bigr)=\begin{cases} \xME/2\pmod1 ,&n\equiv
     4\pmod8\\\xME\quad\pmod1,&n\equiv 0\pmod8\end{cases} 
  \end{equation}
Furthermore, $\cpM\bigl([\cE] \bigr)$ is rational \textnormal{(}lies
in~$\QZ$\textnormal{)}, is independent of metrics and covariant derivatives,
depends only on $[E]\in KO^0(M)$, and is a $\Pp$ bordism invariant.
  \end{conjecture}

\noindent
 Note the shift in the differential pushforward, as in the topological
pushforward~\eqref{eq:22}.  That the $\eta $-invariants satisfy the properties
stated after~\eqref{eq:38} is a theorem.  The key ingredient, already flagged
in~\cite{G}, is that the shift by one implies the characteristic differential
form which computes the variation of the $\eta $-invariant has odd degree,
hence it vanishes.  The other ingredient is the finiteness of the $\Pp$ bordism
group.  (See~\cite{KT2} for more on $\Pp$ bordism.)  \autoref{thm:16} extends
to degrees $n\not\equiv 0\pmod4$, in which case the pushforward is a
topological index.

  \begin{example}[$n=4$]\label{thm:17}
 Let's revisit \autoref{thm:12} with the differential invariant in hand.  The
differential pushforward $\y M\longmapsto \y\cpM(1)$ induces an isomorphism 
  \begin{equation}\label{eq:39}
     \begin{aligned} \Omega ^{\Pp}_4&\longrightarrow \frac{1}{16}\ZZ\!\Bigm/\!\ZZ
      \\ M&\longmapsto \y\xi \mstrut _{M}/2\pmod1\end{aligned} 
  \end{equation}
(See~\cite{St} for more on $\Pp$ $\eta $-invariants.)  In this situation the
differential direct image factorizes through~\eqref{eq:30} to a homomorphism of
topological $KO$-theory groups:
  \begin{equation}\label{eq:40}
     \begin{gathered} \xymatrix@R+.5pc{\cKO^0(M)\ar[r]^{\cp^M}
     \ar@{->>}[d]_<<<<<{\pi  
     _0} & \cKO^{-3}(\pt) \\ KO^0(M)\ar@{-->}[r]^{} &
     KO^{-4}(\pt;\RZ)\ar@{^{(}->}[u]_\iota } 
     \end{gathered} 
  \end{equation}
  \end{example}

In general, the differential pushforward on $\Pp$ manifolds of dimensions
$n\equiv 0,4\pmod 8$ factors as in~\eqref{eq:40} through a topological
homomorphism.  This leads us to ask for topological methods to compute it, and
so to compute the $\eta $-invariant.  A few approaches are described in
\autoref{sec:7}, but we first digress to describe one physics problem in which
these invariants arise.

   \section{Anomalies in M-theory}\label{sec:6}

Very early on, Is foresaw the potential impact of quantum field theory on
mathematics, and he devoted a great deal of effort to foster the interaction
between mathematics and physics.  As part of this effort, he was a progenitor
of the link between index theory and anomalies of spinor fields; the 1984
Atiyah-Singer paper~\cite{AS3} is one of the first to forge this connection.
Soon after, Witten~\cite{W1} brought in $\eta $-invariants.  The conceptual
understanding of anomalies has greatly advanced in the intervening
years~\cite{F2}.  The problem in~\cite{FH2} required computing $\eta
$-invariants on $\Pp$ 12-manifolds of a special sort, as we now explain.
 
M-theory is an 11-dimensional variant of string theory that is closely
connected to 11-dimensional supergravity.  The latter was constructed in~1978
by Cremmer, Julia, and Scherk~\cite{CJS}.  In its Wick-rotated form, there are
four fields: a $\Pp$ structure, a Riemannian metric, a ``Rarita-Schwinger
field'', and a ``$C$-field''.  As explained in~\S\ref{sec:2}, the $\Pp$ field
in the Wick-rotated theory is a consequence of time-reversal symmetry in the
relativistic setting, a crucial property of M-theory.  The Rarita-Schwinger
field is a close cousin of a spinor field, and it gives rise to an anomaly.
There is another anomaly due to the $C$-field, and the main result
of~\cite{FH2} is that these anomalies cancel.  For the present discussion the
$C$-field anomaly is not relevant, but what is important is the global
integrality constraint on the $C$-field.  Namely, the $C$-field is modeled as a
degree four differential cocycle, but with coefficients twisted by the
orientation double cover.  Furthermore, the underlying cohomology class in
$H^4(-;\ZZ_{w_1})$ is constrained~\cite{W2} to have mod~2 reduction equal to
the $4^{\textnormal{th}}$~Stiefel-Whitney class~$w_4$ of the underlying
manifold.  Therefore, the arena for Wick-rotated M-theory is $\Pp$ manifolds
with a $w_1$-twisted integral lift of~$w_4$, topological data we term an
\emph{$\mc$-structure}.
 
Our computational approach to the anomaly question is to compute information
about the bordism group of $\mc$~ 12-manifolds and check anomaly cancellation
on a set of generators.  This is not an off-the-shelf bordism group!  After
much computation and tinkering, we proved this result.

  \begin{theorem}[]\label{thm:18}
  The following six $\mc$-manifolds generate the bordism group $\Omega
^{\mc}_{12}$ localized at the prime~2:  
  \begin{equation}\label{eq:113}
   K\times \HP^{2},\quad
   \RP^{4}\times B, \quad
   (\RP^{4}\#\RP^{4})\times B,\quad 
   W'_{0},\quad
   W''_{0},\quad
   W_{1},
  \end{equation}
where $K$~is the K3 surface, $B$~is an 8-dimensional Bott manifold, and there
are fiber bundles 
  \begin{equation}\label{eq:41}
     \begin{aligned} \HP^2\#\HP^2 \,\longrightarrow\, &W_0'\,\longrightarrow \RP^4
      \\ \RP^8\,\longrightarrow\, &W_0''\longrightarrow \HP^1 \\
      \HP^2\,\longrightarrow \,&W_1\,\longrightarrow \CP^1\times
      \CP^1\end{aligned} 
  \end{equation} 
  \end{theorem}

\noindent
 We refer to~\cite{FH2} for details about these manifolds and about the
M-theory problem.  We computed an appropriate $\eta $-invariant for each of
these six $\Pp$ 12-manifolds, and it is to the computational techniques that we
now turn.

   \section{Topological formul\ae\ for the $\eta $-invariant}\label{sec:7}

For concreteness, specialize to 12-dimensional $\Pp$ manifolds.  Similar
considerations hold in all dimensions $0\pmod 4$ for $\Pp$ manifolds and in
dimensions $2\pmod 4$ for $\Pm$ manifolds.
 
Let $M$~be a closed Riemannian $\Pp$ 12-manifold.  For every real vector bundle
$E\to M$ equipped with an inner product and compatible covariant derivative,
there is an $\eta $-invariant
  \begin{equation}\label{eq:42}
     \xi _M(E)/2\quad \in \QZ
  \end{equation}
According to \autoref{thm:16} it is computed as a pushforward in differential
$KO$-theory, and as in~\eqref{eq:40} the pushforward factors 
  \begin{equation}\label{eq:43}
     \begin{gathered} \xymatrix@R+.5pc{\cKO^0(M)\ar[r]^{\cp^M}
     \ar@{->>}[d]_<<<<<{\pi _0} & \cKO^{-11}(\pt) \\ KO^0(M)\ar@{-->}[r]^{} &
     KO^{-12}(\pt;\RZ)\ar@{^{(}->}[u]_\iota } \end{gathered} 
  \end{equation}
to a homomorphism in topological $KO$-theory.  We indicate three methods to
compute it. 
 
First, if $M$~is oriented and so carries a Spin structure, then 
  \begin{equation}\label{eq:44}
     \xi _M(E)/2 = \ind D_M(E)/2\pmod1, 
  \end{equation}
where $D_M(E)$~is the Spin Dirac operator.\footnote{The usual Dirac operator is
quaternionic, whereas the Clifford linear Dirac operator is real.  The index of
the latter---computed in $KO^{-12}(\pt)$---is what appears in~\eqref{eq:44}; it
is one-half the complex index of the former.}  In other words, in the Spin case
$\xi _M(E)/2$~is a mod~2 index that lifts to an integer index.  As usual, the
integer index has an expression in terms of rational Pontrjagin numbers of~$M$
and $E\to M$. 
 
Second, we indicate a formula of Stolz~\cite[\S5]{St}, derived from a theorem
of Donnelly~\cite{Do}.  Let $\pi \:\hM\to M$ be the orientation double cover,
and let $\sigma \:\hM\to \hM$ be the free orientation-reversing involution.
Suppose $\hM=\partial Z$ for a compact $\Pp$ 13-manifold~$Z$.  Assume $\sigma
$~and $\pi ^*E\to \hM$ both extend over~$Z$ to chosen extensions $\ts$ and
$\sE\to Z$.  If $\ts$~has a finite set~$\{f\}$ of fixed points, then
  \begin{equation}\label{eq:45}
     \xi _M(E)/2 = \sum\limits_{f}\frac{\epsilon _f\tau _f}{2^8}, 
  \end{equation}
where $\tau _f$~is the trace of~$\ts$ lifted to the fiber~$\sE_f$, and
$\epsilon _f=\pm1$ is computed from the action of~$\ts$ on the Pin frames
at~$f$.  
 
These two methods suffice in our problem for each manifold in~\eqref{eq:113}
except~$W_0''$.  For that manifold we used a third method.  It uses a variant
of a theorem of Zhang~\cite{Z}, which he proved analytically based on
techniques in~\cite{BZ}.  To state the theorem, we note that
$\tKO^0(\RP^{20})\cong \zmod{2^{11}}$ with generator $1-[H]$, where $H\to
\RP^{20}$ is the real Hopf line bundle.

  \begin{theorem}[]\label{thm:19}
 Let $M$~be a closed Riemannian $\Pp$ 12-manifold, and suppose $E\to M$ is a
real vector bundle with inner product and orthogonal covariant derivative.
Suppose a smooth map $\y\gamma \:M\longrightarrow \RP^{20}$ satisfies $\gamma
^*w_1(\RP^{20})=w_1(M)$.  Then
  \begin{equation}\label{eq:46}
     \y\gamma _!\bigl([E] \bigr) = 2^{11}\,\frac{\xME}{2}\,\bigl(1-[H] \bigr)
     \qquad \w\textnormal{in }\tKO(\RP^{20}) 
  \end{equation} 
  \end{theorem}

We conclude with a reformulation of \autoref{thm:19} in terms of differential
$KO$-theory.  Namely, the pushforwards
$\y\cpM,\cp^{\RP^{20}}\:\cKO\longrightarrow \cKO(\pt)$ factor through $\pi
_0\:\cKO\rightarrowdbl KO$ in the domain, as in~\eqref{eq:43}.  Using this
factorization and the commutative diagram
  \begin{equation}\label{eq:47}
     \begin{gathered} \xymatrix{\y M\ar[rr]^\gamma \ar[dr]_{\pi ^M} &&
     \y\RP^{20}\ar[dl]^{\pi ^{\RP^{20}}}\\ & \y\pt} \end{gathered} 
  \end{equation}
we reformulate~\eqref{eq:46} as the equality 
  \begin{equation}\label{eq:48}
     \y\cpM = \cp^{\RP^{20}}\circ \gamma _!\:KO^0(M)\longrightarrow
     \cKO^{-11}(\pt)\cong \RZ 
  \end{equation}
This may follow from a more general composition law for pushforwards in
differential theory, perhaps using de Rham currents as in~\cite{FL}, and
perhaps capitalizing as well on the analytic ideas in~\cite{Z}.  It would be
nice to have a more robust theory of generalized differential cocycles in place
to make this---and related---arguments.

 \bigskip\bigskip
\providecommand{\bysame}{\leavevmode\hbox to3em{\hrulefill}\thinspace}
\providecommand{\MR}{\relax\ifhmode\unskip\space\fi MR }
\providecommand{\MRhref}[2]{%
  \href{http://www.ams.org/mathscinet-getitem?mr=#1}{#2}
}
\providecommand{\href}[2]{#2}


\begin{thebibliography}{RSFL10}

\bibitem[AB]{AB}
M.~F. Atiyah and R.~Bott, \emph{A {L}efschetz fixed point formula for elliptic
  complexes. {II}. {A}pplications},
  \href{http://dx.doi.org/10.2307/1970721}{Ann. of Math. (2) \textbf{88}
  (1968)}, 451--491.

\bibitem[ABS]{ABS}
M.~F. Atiyah, R.~Bott, and A.~A. Shapiro, \emph{Clifford modules}, Topology
  \textbf{3} (1964), 3--38.

\bibitem[ADH]{ADH}
Araminta Amabel, Arun Debray, and Peter~J. Haine, \emph{Differential
  Cohomology: Categories, Characteristic Classes, and Connections},
  \href{http://arxiv.org/abs/arXiv:2109.12250}{{\tt arXiv:2109.12250}}.

\bibitem[AmS]{AmS}
W.~Ambrose and I.~M. Singer, \emph{A theorem on holonomy},
  \href{http://dx.doi.org/10.2307/1990721}{Trans. Amer. Math. Soc. \textbf{75}
  (1953)}, 428--443.

\bibitem[APS]{APS}
M.~F. Atiyah, V.~K. Patodi, and I.~M. Singer, \emph{Spectral asymmetry and
  {R}iemannian geometry}, \href{http://dx.doi.org/10.1112/blms/5.2.229}{Bull.
  London Math. Soc. \textbf{5} (1973)}, 229--234.

\bibitem[AS1]{AS1}
M.~F. Atiyah and I.~M. Singer, \emph{The index of elliptic operators. {IV}},
  \href{http://dx.doi.org/10.2307/1970756}{Ann. of Math. (2) \textbf{93}
  (1971)}, 119--138.

\bibitem[AS2]{AS2}
\bysame, \emph{Index theory for skew-adjoint {F}redholm operators}, Inst.
  Hautes {\'E}tudes Sci. Publ. Math. (1969), no.~37, 5--26.

\bibitem[AS3]{AS3}
\bysame, \emph{Dirac operators coupled to vector potentials}, Proceedings of
  the National Academy of Sciences \textbf{81} (1984), no.~8, 2597--2600.

\bibitem[AZ]{AZ}
Alexander Altland and Martin~R. Zirnbauer, \emph{Nonstandard symmetry classes
  in mesoscopic normal-superconducting hybrid structures},
  \href{http://dx.doi.org/10.1103/PhysRevB.55.1142}{Phys. Rev. B \textbf{55}
  (1997)}, 1142--1161.

\bibitem[BNV]{BNV}
Ulrich Bunke, Thomas Nikolaus, and Michael V{\"o}lkl, \emph{Differential
  cohomology theories as sheaves of spectra},
  \href{http://dx.doi.org/10.1007/s40062-014-0092-5}{Journal of Homotopy and
  Related Structures \textbf{11} (2014)}, no.~1, 1--66,
  \href{http://arxiv.org/abs/arXiv:1311.3188}{{\tt arXiv:1311.3188}}.

\bibitem[BS]{BS}
Ulrich Bunke and Thomas Schick, \emph{Smooth {$K$}-theory}, Ast{\'e}risque
  (2009), no.~328, 45--135 (2010).

\bibitem[BSS]{BSS}
Christian Becker, Alexander Schenkel, and Richard~J. Szabo, \emph{Differential
  cohomology and locally covariant quantum field theory},
  \href{http://dx.doi.org/10.1142/S0129055X17500039}{Rev. Math. Phys.
  \textbf{29} (2017)}, no.~1, 1750003, 42,
  \href{http://arxiv.org/abs/arXiv:2011.05768}{{\tt arXiv:2011.05768}}.

\bibitem[BZ]{BZ}
Jean-Michel Bismut and Weiping Zhang, \emph{Real embeddings and eta
  invariants}, \href{http://dx.doi.org/10.1007/BF01444909}{Math. Ann.
  \textbf{295} (1993)}, no.~4, 661--684.

\bibitem[C]{C}
S.~S. Chern, \emph{The geometry of {$G$}-structures},
  \href{http://dx.doi.org/10.1090/S0002-9904-1966-11473-8}{Bull. Amer. Math.
  Soc. \textbf{72} (1966)}, 167--219.

\bibitem[ChS]{ChS}
Jeff Cheeger and James Simons,
  \href{http://dx.doi.org/10.1007/BFb0075216}{\emph{Differential characters and
  geometric invariants}}, Geometry and topology ({C}ollege {P}ark, {M}d.,
  1983/84), Lecture Notes in Math., vol. 1167, Springer, Berlin, 1985,
  pp.~50--80.

\bibitem[CJS]{CJS}
Eugene Cremmer, Bernard Julia, and Joel Scherk, \emph{Supergravity in theory in
  11 dimensions}, Physics Letters B \textbf{76} (1978), no.~4, 409--412.

\bibitem[D]{D}
Freeman~J. Dyson, \emph{The threefold way. {A}lgebraic structure of symmetry
  groups and ensembles in quantum mechanics}, J. Mathematical Phys. \textbf{3}
  (1962), 1199--1215.

\bibitem[De]{De}
Pierre Deligne, \emph{Th{\'e}orie de {H}odge. {II}}, Inst. Hautes {\'E}tudes
  Sci. Publ. Math. (1971), no.~40, 5--57.

\bibitem[Deb]{Deb}
Arun Debray, \emph{Differential cohomology (encyclopedia article)},
  \href{http://arxiv.org/abs/2312.14338}{{\tt 2312.14338}}.

\bibitem[DF]{DF}
Pierre Deligne and Daniel~S. Freed, \emph{Classical field theory}, Quantum
  fields and strings: a course for mathematicians, Vol. 1, 2 (Princeton, NJ,
  1996/1997), Amer. Math. Soc., Providence, RI, 1999, pp.~137--225.

\bibitem[DM]{DM}
Pierre Deligne and John~W. Morgan, \emph{Notes on supersymmetry (following
  {J}oseph {B}ernstein)}, Quantum fields and strings: a course for
  mathematicians, {V}ol. 1, 2 ({P}rinceton, {NJ}, 1996/1997), Amer. Math. Soc.,
  Providence, RI, 1999, pp.~41--97.

\bibitem[Do]{Do}
Harold Donnelly, \emph{Eta invariants for {$G$}-spaces},
  \href{http://dx.doi.org/10.1512/iumj.1978.27.27060}{Indiana Univ. Math. J.
  \textbf{27} (1978)}, no.~6, 889--918.

\bibitem[F1]{F1}
Daniel~S. Freed, \emph{Dirac charge quantization and generalized differential
  cohomology}, Surveys in Differential Geometry, Int. Press, Somerville, MA,
  2000, pp.~129--194.
\href{http://arxiv.org/abs/hep-th/0011220}{{\tt arXiv:hep-th/0011220}}.

\bibitem[F2]{F2}
\bysame, \emph{What is an anomaly?},
  \href{http://arxiv.org/abs/arXiv:2307.08147}{{\tt arXiv:2307.08147}}.

\bibitem[FH1]{FH1}
Daniel~S. Freed and Michael~J. Hopkins, \emph{Reflection positivity and
  invertible topological phases},
  \href{http://dx.doi.org/10.2140/gt.2021.25.1165}{Geom. Topol. \textbf{25}
  (2021)}, no.~3, 1165--1330, \href{http://arxiv.org/abs/arXiv:1604.06527}{{\tt
  arXiv:1604.06527}}.

\bibitem[FH2]{FH2}
\bysame, \emph{{Consistency of M-Theory on Non-Orientable Manifolds}},
  \href{http://dx.doi.org/10.1093/qmath/haab007}{Quart. J. Math. Oxford Ser.
  \textbf{72} (2021)}, no.~1-2, 603--671,
  \href{http://arxiv.org/abs/1908.09916}{{\tt arXiv:1908.09916 [hep-th]}}.

\bibitem[FH3]{FH3}
\bysame, \emph{On Ramond-Ramond fields and $K$-theory}, J. High Energy Phys.
  (2000), \href{http://arxiv.org/abs/hep-th/0002027}{{\tt hep-th/0002027}}.
  Paper 44.

\bibitem[FL]{FL}
Daniel~S. Freed and John Lott, \emph{An index theorem in differential
  {$K$}-theory}, \href{http://dx.doi.org/10.2140/gt.2010.14.903}{Geom. Topol.
  \textbf{14} (2010)}, no.~2, 903--966,
  \href{http://arxiv.org/abs/arXiv:0907.3508}{{\tt arXiv:0907.3508}}.

\bibitem[FM]{FM}
Daniel~S. Freed and Gregory~W. Moore, \emph{Twisted equivariant matter},
  \href{http://dx.doi.org/10.1007/s00023-013-0236-x}{Ann. Henri Poincar{\'e}
  \textbf{14} (2013)}, no.~8, 1927--2023,
  \href{http://arxiv.org/abs/arXiv:1208.5055}{{\tt arXiv:1208.5055}}.

\bibitem[FN]{FN}
Daniel~S. Freed and Andrew Neitzke, \emph{3d spectral networks and classical
  Chern-Simons theory}, Surveys in Differential Geometry \textbf{26} (2021),
  51--155, \href{http://arxiv.org/abs/arXiv:2208.07420}{{\tt
  arXiv:2208.07420}}.

\bibitem[G]{G}
Peter~B. Gilkey, \emph{The eta invariant for even-dimensional {${\rm PIN}_{{\rm
  c}}$} manifolds}, \href{http://dx.doi.org/10.1016/0001-8708(85)90119-7}{Adv.
  in Math. \textbf{58} (1985)}, no.~3, 243--284.

\bibitem[GS]{GS}
Daniel Grady and Hisham Sati, \emph{Differential KO-theory: constructions,
  computations, and applications}, Advances in Mathematics \textbf{384} (2021),
  107671, \href{http://arxiv.org/abs/arXiv:1809.07059}{{\tt arXiv:1809.07059}}.

\bibitem[GY]{GY}
Kiyonori Gomi and Mayuko Yamashita, \emph{Differential {$KO$}-theory via
  gradations and mass terms}, Adv. Theor. Math. Phys. \textbf{27} (2023),
  no.~2, 381--481, \href{http://arxiv.org/abs/arXiv:2111.01377}{{\tt
  arXiv:2111.01377}}.

\bibitem[HHZ]{HHZ}
P.~Heinzner, A.~Huckleberry, and M.R. Zirnbauer, \emph{Symmetry Classes of
  Disordered Fermions},
  \href{http://dx.doi.org/10.1007/s00220-005-1330-9}{Communications in
  Mathematical Physics \textbf{257} (2005)}, no.~3, 725--771,
  \href{http://arxiv.org/abs/arXiv:math-ph/0411040}{{\tt
  arXiv:math-ph/0411040}}.

\bibitem[HS]{HS}
Michael~J. Hopkins and Isadore~M. Singer, \emph{Quadratic Functions in
  Geometry, Topology, and {M}-Theory}, J. Diff. Geom. \textbf{70} (2005),
  329--452,
\href{http://arxiv.org/abs/math/0211216}{{\tt arXiv:math/0211216}}.

\bibitem[K]{K}
Alexei Kitaev, \emph{{Periodic table for topological insulators and
  superconductors}}, \href{http://dx.doi.org/10.1063/1.3149495}{AIP Conf.Proc.
  \textbf{1134} (2009)}, 22--30, \href{http://arxiv.org/abs/0901.2686}{{\tt
  arXiv:0901.2686 [cond-mat.mes-hall]}}.

\bibitem[Kl]{Kl}
K.~R. Klonoff, \emph{An index theorem in differential $K$-theory}, 2008.
  \url{http://repositories.lib.utexas.edu/bitstream/handle/2152/3912/klonoffk16802.pdf?sequence=2}.
  University of Texas Ph.D. thesis.

\bibitem[KT1]{KT1}
R.~C. Kirby and L.~R. Taylor, \emph{Pin structures on low-dimensional
  manifolds}, Geometry of Low-Dimensional Manifolds, 2 ({D}urham, 1989), London
  Math. Soc. Lecture Note Ser., vol. 151, Cambridge Univ. Press, Cambridge,
  1990, pp.~177--242.

\bibitem[KT2]{KT2}
\bysame, \emph{A calculation of {${\rm Pin}^+$} bordism groups}, Commentarii
  Mathematici Helvetici \textbf{65} (1990), no.~1, 434--447.

\bibitem[KZ]{KZ}
Ricardo Kennedy and Martin~R. Zirnbauer, \emph{{Bott Periodicity for
  ${\mathbb{Z}_2}$ Symmetric Ground States of Gapped Free-Fermion Systems}},
  \href{http://dx.doi.org/10.1007/s00220-015-2512-8}{Commun. Math. Phys.
  \textbf{342} (2016)}, no.~3, 909--963,
\href{http://arxiv.org/abs/1409.2537}{{\tt arXiv:1409.2537 [math-ph]}}.

\bibitem[L]{L}
H.~Blaine Lawson, Jr., \emph{{${\rm Spin}^h$} manifolds},
  \href{http://dx.doi.org/10.3842/SIGMA.2023.012}{SIGMA Symmetry Integrability
  Geom. Methods Appl. \textbf{19} (2023)}, Paper No. 012, 7,
  \href{http://arxiv.org/abs/arXiv:2301.09683}{{\tt arXiv:2301.09683}}.

\bibitem[LM]{LM}
H.~B. Lawson, Jr. and M.-L. Michelsohn, \emph{Spin geometry}, Princeton
  Mathematical Series, vol.~38, Princeton University Press, Princeton, NJ,
  1989.

\bibitem[S]{S}
I.~M. Singer, 2010. \url{https://www.youtube.com/watch?v=5FoaMcCJnmQ}.
  Interview of Isadore Singer for the MIT+150 Infinite History Project.

\bibitem[SRFL]{SRFL}
Shinsei Ryu, Andreas~P. Schnyder, Akira Furusaki, and Andreas~W.W. Ludwig,
  \emph{{Topological insulators and superconductors: Tenfold way and
  dimensional hierarchy}},
  \href{http://dx.doi.org/10.1088/1367-2630/12/6/065010}{arXiv:0912.2157; New
  J.Phys. \textbf{12} (2010)}, 065010,
\href{http://arxiv.org/abs/arXiv:0912.2157}{{\tt arXiv:0912.2157}}.

\bibitem[SS1]{SS1}
James Simons and Dennis Sullivan, \emph{Axiomatic characterization of ordinary
  differential cohomology}, \href{http://dx.doi.org/10.1112/jtopol/jtm006}{J.
  Topol. \textbf{1} (2008)}, no.~1, 45--56,
  \href{http://arxiv.org/abs/arXiv:math/0701077}{{\tt arXiv:math/0701077}}.

\bibitem[SS2]{SS2}
\bysame, \emph{Structured vector bundles define differential {$K$}-theory},
  Quanta of maths, Clay Math. Proc., vol.~11, Amer. Math. Soc., Providence, RI,
  2010, pp.~579--599. \href{http://arxiv.org/abs/arXiv:0810.4935}{{\tt
  arXiv:0810.4935}}.

\bibitem[St]{St}
Stephan Stolz, \emph{Exotic structures on 4-manifolds detected by spectral
  invariants}, Inventiones mathematicae \textbf{94} (1988), no.~1, 147--162.

\bibitem[SW]{SW}
R.~F. Streater and A.~S. Wightman, \emph{P{CT}, spin and statistics, and all
  that}, Princeton Landmarks in Physics, Princeton University Press, Princeton,
  NJ, 2000. Corrected third printing of the 1978 edition.

\bibitem[W1]{W1}
Edward Witten, \emph{Global gravitational anomalies},
  \href{http://dx.doi.org/10.1007/BF01212448}{Commun. Math. Phys. \textbf{100}
  (1985)},
197.

\bibitem[W2]{W2}
\bysame, \emph{On flux quantization in {$M$}-theory and the effective action},
  \href{http://dx.doi.org/10.1016/S0393-0440(96)00042-3}{J. Geom. Phys.
  \textbf{22} (1997)}, no.~1, 1--13,
  \href{http://arxiv.org/abs/arXiv:hep-th/9609122}{{\tt arXiv:hep-th/9609122}}.

\bibitem[WS]{WS}
Chong Wang and T~Senthil, \emph{Interacting fermionic topological
  insulators/superconductors in 3D}, Physical Review B \textbf{89} (2014),
  no.~19, 195124, \href{http://arxiv.org/abs/arXiv:1401.1142}{{\tt
  arXiv:1401.1142}}.

\bibitem[Z]{Z}
WeiPing Zhang, \emph{A {$\rm{mod}\,2$} index theorem for {$\rm pin^-$}
  manifolds}, \href{http://dx.doi.org/10.1007/s11425-016-9040-7}{Sci. China
  Math. \textbf{60} (2017)}, no.~9, 1615--1632,
  \href{http://arxiv.org/abs/arXiv:1508.02619}{{\tt arXiv:1508.02619}}.

\end{thebibliography}
  \end{document}